\documentclass[a4paper,11pt]{amsart}

\newtheorem{theorem}{Theorem}[section]
\newtheorem{lemma}[theorem]{Lemma}
\newtheorem{proposition}[theorem]{Proposition}

\newcommand{\ol}[1]{\overline{#1}}
\newcommand{\ul}[1]{\underline{#1}}
\newcommand{\GreenL}{\mathrel{\mathcal{L}}}
\newcommand{\GreenR}{\mathrel{\mathcal{R}}}
\newcommand{\GreenH}{\mathrel{\mathcal{H}}}

\newcommand{\GreenJ}{\mathrel{\mathcal{J}}}

\newcommand{\BZ}{\mathbb{Z}}

\begin{document}

    \title{Wreath product decompositions for \\ triangular matrix semigroups}

    \author{MARK KAMBITES}

    \address{Fachbereich Mathematik / Informatik, Universit\"at Kassel,
             34109 Kassel,  Germany
    }
    \email{kambites@theory.informatik.uni-kassel.de}

    \author{BENJAMIN STEINBERG}

    \address{School of Mathematics and Statistics, Carleton University,
             Ottawa, Ontario, K1S 5B6, Canada
    }
    \email{bsteinbg@math.carleton.ca}

    \maketitle
\begin{abstract}
We consider wreath product decompositions for semigroups of
triangular matrices. We exhibit an
explicit wreath product decomposition for the semigroup of all
$n \times n$ upper triangular matrices over a given
field $k$, in terms of aperiodic semigroups and affine groups over $k$. In the
case that $k$ is finite this decomposition is optimal, in the sense
that the number of group terms is equal to the group complexity of
the semigroup. We also obtain some decompositions for semigroups of
triangular matrices over more general rings and semirings.
\end{abstract}

\section{Introduction}\label{sec_intro}

Some of the most natural and frequently occurring semigroups are
those of upper triangular matrices over a given ring or field. For
example, such semigroups arise in the study of algebraic semigroups,
where Putcha \cite{Putcha88} has proven that a connected algebraic
monoid with zero over a field has a faithful rational triangular
representation if and only if its group of units is
solvable \cite{Putcha88}. It follows that triangularizable monoids
can be thought of as a natural generalisation of solvable groups.
More recently, Almeida, Margolis and Volkov \cite{Almeida05} have
shown that semigroups of triangular matrices over finite fields
generate natural pseudovarieties. Almeida, Margolis, Steinberg and
Volkov \cite{words,radical} have since considered arbitrary fields
and have obtained language-theoretic consequences. Further
properties of these semigroups have been described by
Okninski \cite{Okninski98}.

Perhaps the most productive approach to the study of finite
semigroups is through coverings by wreath products. In the 1960s,
Krohn and Rhodes \cite{KR1,Krohn68,Arbib68} showed that every finite
semigroup can be expressed as a \textit{divisor} (a homomorphic
image of a subsemigroup) of a wreath product of finite groups and
finite aperiodic monoids. The \textit{group complexity} of a finite
semigroup is the smallest number of group terms in such a
decomposition, and is a key concept in finite semigroup theory.

In a previous article \cite{Kambites05triangular}, the first author
computed the group complexity of the semigroup $T_n(k)$ of all $n
\times n$ upper triangular matrices over a given finite field $k$,
and of certain related semigroups. However, the methods used did not
result in explicit wreath product decompositions. The main objective of this
article is to establish an explicit wreath product decomposition for
each semigroup of the form $T_n(k)$, and hence for every semigroup
of triangular matrices over a finite field. This decomposition is
optimal, in the sense that the number of group terms in the
decomposition is equal to the group complexity of $T_n(k)$. Moreover
every group appearing is a product of subgroups of $T_n(k)$.

In the process, we obtain some results applicable in a more general context.
While Krohn-Rhodes theory is traditionally concerned with finite semigroups,
there have been numerous attempts to extend it to well-behaved classes of
infinite semigroups \cite{Birget84,Elston02,Henckell88}. Our method for decomposing $T_n(k)$
is fully applicable in the case that the field $k$, and hence also the
semigroup $T_n(k)$, is infinite. We also obtain some wreath product
decompositions, although not in terms of groups and aperiodic semigroups, for
triangular matrix semigroups over more general rings and semirings with identity.

In addition to this introduction, this paper comprises four sections.
In Section~\ref{sec_krohnrhodes}, we briefly recall the key definitions
and results of Krohn-Rhodes theory as applied to abstract monoids, including
division, wreath products, the Prime Decomposition Theorem
and group complexity. Section~\ref{sec_triangular} introduces triangular
matrix semigroups, and briefly describes their structure, before
reviewing the results of the first
author \cite{Kambites05triangular} characterising their group complexity.

Section~\ref{sec_decomposition} contains the main original results
of the paper; we obtain an explicit decomposition for each semigroup
$T_n(k)$, and hence for every semigroup of triangular matrices over
a field, as a wreath product of aperiodic monoids and affine groups
over $k$. We also obtain some related decompositions for triangular
matrix semigroups over rings and semirings with identity. Finally,
in Section~\ref{sec_comparison}, we compare our results with those
which can be obtained using a standard decomposition method of
Eilenberg and Tilson \cite{TilsonXI}; the latter produces a
suboptimal decomposition for $T_n(k)$, but alternative optimal
decompositions for certain important divisors.

Throughout this paper,  all functions are applied on the right of
their arguments. If $S$ and $T$ are sets then we denote by $S^T$ the
set of all functions from $T$ to $S$.  We assume familiarity with
the standard terminology, notation and foundational results of
structural semigroup theory; a detailed introduction to these is
given by Howie \cite{Howie95}. By contrast, we assume no prior
knowledge which is particular to the study of finite semigroups; we
intend that this article should be fully accessible to the reader
with experience only of infinite semigroups.

\section{Wreath Products, Division and Complexity}\label{sec_krohnrhodes}

In this section, we briefly introduce the basic concepts of wreath
products, division and complexity. We restrict ourselves to the special
case of abstract monoids (as opposed to transformation semigroups),
since this suffices for our purpose. A detailed and more general
introduction is given by Eilenberg \cite{Eilenberg76}.

Let $S$ and $T$ be semigroups. We say that $S$ \textit{divides} $T$,
and write $S \prec T$, if $S$ is a homomorphic image of some
subsemigroup of $T$. The relation of division is easily verified to
be reflexive and transitive.

Let $S$ and $T$ be monoids.  Then $S^T$ is a monoid with pointwise
product:  if $f,g\in S^T$ and $t\in T$, then by definition $t(fg)=
(tf)(tg)$. There is also a natural left action of $T$ on $S^T$
defined as follows: if $f\in S^T$, $t_1,t_2\in T$, then
${}^{t_1}\!{f}:T\to S$ is given by $$t_2{}^{t_1}\!{f}=(t_2t_1)f.$$
Then the \textit{wreath product} of $S$ and $T$, denoted $S\wr T$,
is the monoid with underlying set $S^T \times T$, and multiplication
given by
$$(f, a) (g, b) = (f{}^{a}\!{g}, ab).$$
The wreath product of monoids is not associative; however, $(S_3\wr S_2)\wr S_1$
is isomorphic to a submonoid of $S_3\wr (S_2\wr S_1)$.  For this reason,
we define the \textit{iterated wreath product} of a sequence of three or more
monoids inductively by
$$S_n \wr S_{n-1} \wr \cdots \wr S_1 \ = \ S_n \wr \left( S_{n-1} \wr \cdots \wr S_1 \right)$$
so as to obtain the largest monoid possible.

Recall that a semigroup is called \textit{aperiodic} if it has no
non-trivial subgroups. In the following proposition we state without
proof a few well-known properties of the wreath product which we
shall need.
\begin{proposition}\label{prop_wreathproperties}
Let $A$, $B$, $C$ and $D$ be finite monoids.
\begin{itemize}
\item[(i)] If $A \prec B$ then $A \wr C \prec B \wr C$ and $C \wr A \prec C \wr B$.
\item[(ii)] $A \times B \prec A \wr B$.
\item[(iii)] $(A \wr B) \times (C \wr D) \prec (A \times C) \wr (B \times D)$.
\item[(iv)] If $A$ and $B$ are groups then $A \wr B$ is a group.
\item[(v)] If $A$ and $B$ are aperiodic then $A \wr B$ is aperiodic.
\end{itemize}
\end{proposition}

We shall also need an elementary decomposition that is perhaps not
so well known:

\begin{proposition}\label{prop absorbprod}
Let $A$, $B$, $C$ be monoids.  Then $(A\wr B)\times C$ embeds in
$A\wr (B\times C)$.
\end{proposition}
\begin{proof}
First we define a homomorphism $\alpha:A^B\to A^{B\times C}$ by
$(b,c)f\alpha = bf$.   Next we define $\psi:(A\wr B)\times
C\rightarrow A\wr (B\times C)$ by $$((f,b),c)\psi =
(f\alpha,(b,c)).$$  Let us verify that $\psi$ is a homomorphism.
\begin{align*}
((f,b),c)\psi \ ((g,b'),c')\psi &= (f\alpha,(b,c)) \ (g\alpha,(b',c')) \\
& = (f\alpha{}^{(b,c)}\!{g\alpha},(bb',cc')).
\end{align*}
But for $(b_0,c_0)\in B\times C$,
$$(b_0,c_0){}^{(b,c)}\!{g\alpha} = (b_0b,c_0c)g\alpha = b_0bg = b_0{}^b\!{g}.$$
Thus ${}^{(b,c)}\!{(g\alpha)} = ({}^b\!{g})\alpha$.  Hence we may
conclude
\begin{align*}
((f,b),c)\psi \ ((g,b'),c')\psi  &=
(f\alpha({}^b\!{g})\alpha,(bb',cc'))\\ &= ((f{}^{b}\!{g})\alpha,(bb',cc'))\\
&= [((f,b),c)((g,b'),c')]\psi
\end{align*}
 and so $\psi$ is a
homomorphism.  It is clear that $\psi$ is injective.
\end{proof}

 Let $X$ be a finite set. We denote by
$\widetilde {X}$ the monoid consisting of the identity map and all
constant maps on $X$; clearly, $\widetilde {X}$ is an aperiodic
monoid. Now if $A$ is a monoid of transformations of $X$, then the
\textit{augmented monoid} $\overline{A}$ of $A$ with respect to its
action on $X$ is the monoid generated by transformations in $A$ and
those in $\widetilde {X}$. The following proposition, a
proof of which can be found in Eilenberg \cite{Eilenberg76}, provides a
decomposition of an augmented monoid in terms an
aperiodic monoid and the underlying monoid.

\begin{proposition}\label{prop_augmentationdecomp}
Let $A$ be a finite monoid of transformations of a set $X$. Then
$\overline {A}\prec \widetilde {X} \wr A$.
\end{proposition}

The importance of wreath products for the study of finite semigroups
stems from the following structure theorem of Krohn and
Rhodes \cite{KR1,Arbib68}.
\begin{theorem}\textsf{(The Prime Decomposition Theorem, Krohn-Rhodes 1968)}\ Let $S$ be a finite semigroup. Then $S$ divides some
iterated wreath product each of whose terms is either (i) a finite
simple group which divides $S$ or (ii) a finite aperiodic monoid.
\end{theorem}
A \textit{Krohn-Rhodes decomposition} for a semigroup $S$ is an expression of $S$
as a divisor of an iterated wreath product of groups and aperiodic monoids.
Given such a decomposition for $S$,
Proposition~\ref{prop_wreathproperties} tells us that we can
combine adjacent groups terms and adjacent aperiodic terms to obtain
an \textit{alternating} decomposition of the form:
$$A_n \wr G_n \wr A_{n-1} \wr \dots \wr A_1 \wr G_n \wr A_0$$
where each $A_i$ is aperiodic, each $G_i$ is a group, and all terms
except possibly $A_0$ and $A_n$ are non-trivial. (Note, though, that in
doing so we may lose the property that the group terms are divisors of
$S$.) The number $n$, that
is, the number of group terms, is called the
\textit{group length} of the decomposition. A natural structural constant
which can be associated with a finite semigroup
$S$ is the minimal group length of a Krohn-Rhodes decomposition for $S$;
this number is
called the \textit{group complexity} of $S$. A
decomposition for $S$ is said to be \textit{optimal} if its group length
equals the group complexity of $S$.

Much effort has been put into the study of Krohn-Rhodes decompositions,
and in particular of certain algorithmic problems. Various algorithms have been
developed for finding wreath product decompositions for semigroups; some of
these will be discussed in Section~\ref{sec_comparison} below. A major
open question is that of whether group complexity is \textit{decidable},
that is, whether there is an algorithm which, given the multiplication table
for a finite semigroup $S$, determines the group complexity of $S$.

We remark briefly upon the relationship between these two problems,
and in particular on the implications of the latter for the former.
In theory, knowing the group complexity of a finite semigroup allows
one to compute an optimal decomposition. Indeed, if ones knows that
a semigroup $S$ admits a decomposition of group length $n$, then one can
in principle enumerate multiplication tables of divisors of
alternating wreath products with $n$ group terms, and test them for
isomorphism with $S$. In practice, of course, this algorithm is
completely infeasible -- the cardinality of an iterated
wreath product grows extremely fast as function of
the cardinalities of the terms, and no sensible upper bounds are
known even on the latter. Hence, situations can arise in which the
complexity of a semigroup is known, but an explicit optimal decomposition is
not.
Indeed, the following key result
of Rhodes \cite{Rhodes74} often gives rise to such situations.
\begin{theorem} \textsf{(The Fundamental Lemma of Complexity, Rhodes
1974)}\
Let $S$ and $T$ be finite semigroups, and suppose there exists a surjective
morphism $S \to T$ which is injective when restricted to each subgroup of
$S$. Then $S$ and $T$ have the same group complexity.
\end{theorem}
The Fundamental Lemma is an extremely powerful tool for computing
the group complexity of a semigroup. The proof of the Lemma given by
Tilson \cite{TilsonXII} is constructive in the sense that, given an
optimal wreath product decomposition for a semigroup $T$ and a
surjective morphism $S \to T$ which is injective on subgroups, it
does provide an optimal decomposition for $S$.  However the
construction is quite involved and in practice it is hard to see what
groups and aperiodic monoids appear.

\section{Triangular Matrix Semigroups}\label{sec_triangular}

Let $R$ be a semiring with identity $1$ and zero $0$.
If $x$ is an $n \times n$ matrix then for $1 \leq i, j \leq n$ we
denote by $x_{ij}$ the entry of $x$ in position $(i,j)$, that
is, in the $i$th row and $j$th column, of $x$. Recall that the matrix $x$
is \textit{(upper) triangular} if $x_{ij} = 0$ whenever $1 \leq j < i \leq
n$. We call an upper triangular matrix \textit{(upper) unitriangular} if,
in addition, $x_{ii} = 0$ or $x_{ii} = 1$ for $1 \leq i \leq n$. We call
$x$ a \textit{subidentity} if it is unitriangular and $x_{ij} = 0$ whenever
$i \neq j$. We denote by $T_n(R)$ and $UT_n(R)$ the semigroups of all
$n \times n$ upper triangular matrices and of all $n \times n$ unitriangular
matrices respectively, with entries drawn from $R$, the operation in both
cases being usual matrix multiplication. Note that $T_1(R)$ is just the
multiplicative semigroup of $R$.

We shall be especially interested in the case that the semiring $R$ is a field
$k$. In this case, we define a relation $\sigma$ on each semigroup $T_n(k)$ by $x \mathrel{\sigma} y$ if and only $x = \lambda y$ for some
non-zero scalar $\lambda$. This relation is easily verified to be a congruence on
$T_n(k)$. The \textit{projective triangular semigroup} $PT_n(k)$ is the
quotient semigroup $T_n(k) / \sigma$; we denote by $\overline{x}$ the element
of $PT_n(k)$ which is the $\sigma$-equivalence class of a matrix $x \in T_n(k)$.

The group of units of $T_n(k)$ [respectively, $UT_n(k)$, $PT_n(k)$]
is denoted $T_n^*(k)$ [$UT_n^*(k)$, $PT_n^*(k)$]. It consists of
those triangular matrices whose diagonal entries are non-zero
[respectively, triangular matrices whose diagonal entries are $1$,
equivalence classes of triangular matrices whose diagonal entries
are non-zero]. Note that $T_1^*(k)$ is the multiplicative group
of the field $k$.

We introduce a notion of upper triangular row and column operations on
$T_n(R)$ and hence on
$UT_n(R)$. By a \textit{row operation} on an upper triangular matrix
we shall mean either (i) adding a multiple of one row to a row
\textit{above} or (ii) scaling a row by an element of $R$. There
is an obvious analogous definition of \textit{column operations} of
different types, a type (i) operation being adding a multiple of one
column to a column to the \textit{right}. The following easy proposition
characterises Green's relations in $T_n(R)$ and $UT_n(R)$ in terms of these
operations.
\begin{proposition}\label{prop_greens_characterization}
Let $n$ be a positive integer and $R$ a semiring with identity. Two matrices in
$T_n(R)$ [respectively, $UT_n(R)$] are:
\begin{itemize}
\item[(i)] $\mathcal{L}$-related exactly if each can be obtained from the other
 by [unitriangular] row operations;
\item[(ii)] $\mathcal{R}$-related exactly if each can be obtained from the
 other by [unitriangular] column operations;
\item[(iii)] $\mathcal{J}$-related exactly if each can be obtained from the
 other by [unitriangular] row and column operations.
\end{itemize}
\end{proposition}

We now turn our attention to the case of a finite field $k$. The
following proposition, parts of which go back at least as far as
Putcha \cite{Putcha88}, characterises the regular elements in
$T_n(k)$. A proof can be found in a previous article
of the first author \cite{Kambites05triangular}.

\begin{proposition}\label{prop_regular_characterization}
Let $n$ be a positive integer and $k$ a finite field.
Let $x \in T_n(k)$ or $x \in UT_n(k)$. Then the following are equivalent:
\begin{itemize}
\item[(i)] $x$ is regular;
\item[(ii)] every row in $x$ is a linear combination of rows in $x$ with
 non-zero diagonal entries;
\item[(iii)] every column in $x$ is a linear combination of columns in $x$
 with non-zero diagonal entries;
\item[(iv)] $x$ is $\mathcal{J}$-related to a subidentity.
\end{itemize}
\end{proposition}

Factoring out a monoid by a subgroup of the group of units that is
central in the monoid gives rise to a congruence contained in
$\GreenH$. The following simple observation is a special case of
well-known and elementary facts about congruences contained in
$\GreenH$.
\begin{proposition}\label{prop_ptstructure}
Let $n$ be a positive integer, $k$ a field and $x,y \in T_n(k)$. Then
\begin{itemize}
\item[(i)] $x$ is regular in $T_n(k)$ if and only if $\ol{x}$ is regular in $PT_n(k)$;
\item[(ii)] $x \GreenL y$ in $T_n(k)$ if and only if $\ol{x} \GreenL \ol{y}$ in $PT_n(k)$;
\item[(iii)] $x \GreenR y$ in $T_n(k)$ if and only if $\ol{x} \GreenR \ol{y}$ in $PT_n(k)$;
\item[(iv)] $x \GreenJ y$ in $T_n(k)$ if and only if $\ol{x} \GreenJ \ol{y}$ in $PT_n(k)$.
\end{itemize}
\end{proposition}

We recall the following theorem of the first author \cite{Kambites05triangular}.
\begin{theorem}\label{thm_trimatcomplexity}\textsf{(Kambites 2004)}\
Let $n$ be a positive integer, and $k$ a finite field. If $n > 1$ or
$k = \BZ_2$ then $T_n(k)$, $UT_n(k)$ and $PT_n(k)$ have complexity $n-1$.
If $n = 1$ and $k \neq \BZ_2$ then $UT_n(k)$ and $PT_n(k)$ have complexity
$0$, while $T_n(k)$ has complexity $1$.
\end{theorem}
We remark that the scope of this result has since been extended by
Mintz \cite{Mintz05}; he observes that triangular matrix semigroups form
a special class of \textit{quiver algebra} and that the result extends
naturally to cover a somewhat larger class of quiver algebras.

The proof of Theorem~\ref{thm_trimatcomplexity} is somewhat
technical, and makes extensive use of the Fundamental Lemma of
Complexity, both directly and through the application of a result of
Rhodes and Tilson \cite{Rhodes69}. Consequently, it
does not give rise to explicit Krohn-Rhodes decompositions
for the semigroups in question. In the next section, we shall show
how to obtain such decompositions for semigroups of the form
$T_n(k)$, and hence for every triangular matrix semigroup over a
field. 

\section{Decompositions for Triangular Matrix Semigroups}\label{sec_decomposition}

Our main objective in this section is to compute an explicit
decomposition
for each semigroup of the form $T_n(k)$ with $k$ a field, as a divisor
of an alternating wreath product of groups and aperiodic monoids. In the case
that $k$ is finite, this decomposition will be optimal, in the sense that
its group length equals the group complexity of the semigroup as described
by Theorem~\ref{thm_trimatcomplexity}. In the process, we also obtain some
decompositions for triangular matrix semigroups over more general rings
and semirings.

Let $R$ be a semiring and $n$ a positive integer. We consider the
$R$-module $R^n$ of $1\times n$ row vectors over $R$. Recall that an
\textit{affine transformation} of $R^n$ is a map of the form $v
\mapsto vX + c$ for some $n \times n$ matrix $X$ and some vector $c
\in R^n$. We say that the transformation is \textit{affine (upper)
triangular} if $X$ is upper triangular, and \textit{affine scaling}
if $X$ is of the form $\lambda I$ where $\lambda \in R$ and $I$ is
the identity matrix.

The \textit{affine monoid} $A_n(R)$ of degree $n$ over $R$ is the
monoid of all affine transformations of $R^n$, with operation
composition. It is readily verified that the sets of affine
triangular and affine scaling maps form submonoids; these we call
the \textit{affine triangular monoid} $AT_n(R)$ and the
\textit{affine scaling monoid} $AS_n(R)$ respectively. The
\textit{affine group} $A_n^*(R)$, the \textit{affine triangular
group} $AT_n^*(R)$ and the \textit{affine scaling group} $AS_n^*(R)$
are the groups of units of $A_n(R)$, $AT_n(R)$ and $AS_n(R)$
respectively.  We remark that the various affine groups are
semidirect products of the appropriate matrix groups and with the
additive group of translations.

There is a natural embedding of an affine triangular monoid of degree
$n-1$ into an upper triangular monoid of degree $n$.
\begin{proposition}\label{prop embedding}
Let $n \geq 2$ and let $R$ be a semiring. Then $AT_{n-1}(R)$ and
$AS_{n-1}(R)$ embed in $T_n(R)$.
\end{proposition}
\begin{proof}
From the definition, $AS_{n-1}(R)$ is a subsemigroup of
$AT_{n-1}(R)$, so it suffices to show that the latter embeds in
$T_n(R)$. Given an affine triangular map $f$ given by $v \mapsto vX
 + c$ we define an $n \times n$ matrix
$$M_f = \left( \begin{matrix} 1 & c \\ 0 & X \end{matrix} \right).$$
That the matrix $M_f$ is upper triangular follows from the fact that
$X$ is upper triangular. If we identify $v\in R^{n-1}$ with $(1,v)$
then it is routine to verify that $(1,v)M_f = (a,vf)$ and so
$f\mapsto M_f$ gives an embedding of $AT_{n-1}(R)$ into $T_n(R)$, as
required.
\end{proof}

The following lemma is the main inductive step in our
decompositions.  If $X$ is a matrix, we write $X^{T}$ for its
transpose.
\begin{lemma}\label{lemma_induction}
Let $n \geq 2$ and $R$ be a semiring with identity. Then
$$T_n(R) \prec \left[ AS_{n-1}(R) \wr T_{n-1}(R) \right] \times T_1(R).$$
\end{lemma}
\begin{proof}
We view each $s \in T_n(R)$ as a block matrix
$$s = \left( \begin{matrix} M_s & v_s \\ 0 & c_s \end{matrix} \right)$$
where $M_s$ is an $(n-1) \times (n-1)$ matrix which clearly lies in
$T_{n-1}(R)$, $v_s$ is an $n \times 1$ column vector and $c_s$ is a
$1 \times 1$ matrix. Now we define
$$\psi : T_n(R) \to [AS_{n-1}(R) \wr T_{n-1}(R)] \times T_1(R)$$
by
$$s \psi = (f_s,  M_s, c_s)$$
where for every $M \in T_{n-1}(R)$, the element $M f_s \in
AS_{n-1}(R)$ is given by
$$w (X f_s) = (X v_s + w^{T} c_s)^{T}.$$
Clearly, $\psi$ is well-defined; it is also injective, since for any
$s$, we have \mbox{$v_s = [\ul{0} (If_s)]^{T}$} where $I \in
T_{n-1}(R)$ is the identity matrix and $\ul{0} \in R^{n-1}$ is the
zero vector.

To prove the lemma, it will now suffice to show that $\psi$ is a
homomorphism. Since $$\left(\begin{matrix} M_s& v_s \\ 0 &
c_s\end{matrix}\right)\left(\begin{matrix} M_t& v_t \\ 0 &
c_t\end{matrix}\right)=\left(\begin{matrix} M_sM_t& M_sv_s+v_sc_t \\
0 & c_sc_t\end{matrix}\right),$$ we have: $M_{st} = M_s M_t$,
$c_{st} = c_s c_t$ and $v_{st} = M_sv_t + v_sc_t$.   So, recalling
the definition of the wreath product, it remains to show that
$f_{st} = f_s{}^{M_s}\!{f_t}$.  That is we must show \mbox{$w(X
f_{st}) = w [(X f_s) (X M_s f_t)]$} for all $w \in R^n$ and $X \in
T_{n-1}(R)$. But
\begin{align*}
w (X f_{st}) &= (X v_{st} + w^{T} c_{st})^{T} \\
             &= (X (M_s v_t + v_s c_t) + w^{T} c_s c_t)^{T} \\
             &= (X M_s v_t + X v_s c_t + w^{T} c_s c_t)^{T} \\
             &= (X M_s v_t + (X v_s + w^{T} c_s) c_t)^{T} \\
             & =(X M_s v_t)^{T} + (Xv_s +w^{T}c_s)^{T}c_t \\
             & = (X M_s v_t)^{T}  + (wXf_s)c_t \\
             &=  (XM_sv_t + (w(Xf_s))^{T}c_t)^{T}\\
             &= w [(X f_s) (X M_s f_t)]
\end{align*}
as required.
\end{proof}

Lemma~\ref{lemma_induction} leads easily to the following decomposition
for $T_n(R)$ in terms of affine scaling monoids and the multiplicative
semigroup of $R$.
\begin{theorem}\label{thm_ringdecomp}
Let $n \geq 2$ and $R$ be a semiring with identity. Then $$T_n(R)
\prec  AS_{n-1}(R) \wr AS_{n-2}(R) \wr \dots \wr (AS_1(R) \times
T_1(R)^n).$$
\end{theorem}
\begin{proof}
We use induction on $n$. When $n=2$ then using
Lemma~\ref{lemma_induction} and Proposition \ref{prop absorbprod} we
have
$$T_2(R) \prec  \left[ AS_1(R) \wr T_1(R) \right] \times T_1(R)\prec AS_1(R) \wr T_1(R)^2$$
as required. Now let $n \geq 3$ and assume true for smaller $n$.
Then again using Lemma~\ref{lemma_induction} and Proposition
\ref{prop absorbprod} we obtain
\begin{align*}
T_n(R) &\prec  \left[ AS_{n-1}(R) \wr T_{n-1}(R) \right] \times T_1(R) \\
       &\prec  \left[ AS_{n-1}(R) \wr \left(AS_{n-2}(R) \wr \dots \wr (AS_1(R) \times T_1(R)^{n-1})\right) \right] \times T_1(R) \\
       &\prec  AS_{n-1}(R) \wr AS_{n-2}(R) \wr \dots \wr (AS_1(R)
       \times
       T_1(R)^n)
\end{align*}
as required.
\end{proof}

As a consequence of Theorem~\ref{thm_ringdecomp}, we obtain a group
length $n-1$ decomposition for each semigroup $T_n(k)$ with $k$ a
field.
\begin{theorem}\label{thm_fielddecomp}
Let $n \geq 2$ and $k$ be a field. Then $T_n(k)$ divides
$$\widetilde{k^{n-1}} \wr AS^*_{n-1}(k) \wr \widetilde{k^{n-2}} \wr AS^*_{n-2}(k) \wr \dots \wr \widetilde{k} \wr \left[ AS^*_1(k) \times T_1^*(k)^n \right] \wr U_1^n$$
where $U_1$ is the two-element semilattice.
\end{theorem}
\begin{proof}
By Theorem~\ref{thm_ringdecomp} we have that
$$T_n(k) \prec AS_{n-1}(k) \wr AS_{n-2}(k) \wr \dots \wr (AS_1(k) \times T_1(k)^n).$$
For each $i$, it is easily seen that the affine monoid $AS_i(k)$
consists precisely of $AS^*_i(k)$ and constant maps on $k^i$; hence,
$AS_i(k)$ is the augmented monoid of $AS^*_i(k)$ with respect to its
action on $k^i$, and so by Proposition~\ref{prop_augmentationdecomp}
we have
$$AS_i(k) \prec \widetilde{k^i}  \wr AS^*_i(k).$$
Also, it is easy to see that the group with zero $T_1(k)$ divides
$T_1^*(k) \times U_1$. It follows that
$$AS^*_1(k)\times T_1(k)^n \prec AS^*_1(k)\times T_1^*(k)^n\times U_1^n\prec  (AS^*_1(k)\times T_1^*(k)^n)\wr U_1^n.$$
The result is now clear.
\end{proof}

Recall that the \textit{pseudovariety} generated by a finite semigroup $S$
is the class of all divisors of finite direct products of $S$. In general,
a finite semigroup $S$ does not necessarily admit an optimal Krohn-Rhodes decomposition
whose group terms are divisors of $S$, or even in the pseudovariety generated
$S$. Here
we have succeeded in finding for $T_n(k)$ an optimal Krohn-Rhodes
decomposition in which every group is a subgroup of the group of
units $T_n(k)$ except one, which is a direct product of two
subgroups of $T^*_n(k)$. Indeed, Proposition \ref{prop embedding}
implies that each $AS^*_{m}(k)$ with $1\leq m\leq n-1$ embeds in
$T^*_n(k)$. On the other hand $T_1^*(k)^n$ is just the
diagonal subgroup of $T_n^*(k)$.

\section{Comparison with Depth Decomposition}\label{sec_comparison}

Considerable thought has been put into algorithmic methods for
obtaining explicit Krohn-Rhodes decompositions for finite
transformation semigroups. The original proof of Krohn and
Rhodes \cite{Krohn68} is essentially algorithmic; however, the
decompositions it yields are far from optimal. A substantial
improvement is the \textit{holonomy method}, which was developed by
Eilenberg \cite{Eilenberg76}, in conjunction with Tilson, using
techniques of Zeiger \cite{Zeiger67} and Ginzburg \cite{Ginzburg68};
see also Holcombe \cite{Holcombe81} for a good exposition with a
small correction to Eilenberg's definitions.

When attention is restricted to abstract semigroups (as opposed to
transformation semigroups), better methods are available. The
\textit{depth decomposition} method of Eilenberg and
Tilson \cite{TilsonXI} is known to yield decompositions for abstract
semigroups which are at least as short as, and sometimes shorter
than, holonomy decompositions. We briefly recall the depth
decomposition method; for full details, see Tilson \cite{TilsonXI}.

Recall that a $\mathcal{J}$-class is called \textit{essential} if it contains
a non-trivial subgroup. The \textit{depth} of an essential
$\mathcal{J}$-class is the length of the longest chain of essential
$\mathcal{J}$-classes strictly above it. The \textit{depth of the semigroup} is defined
to be the length of the longest chain of essential $\mathcal{J}$-classes
in the semigroup, that is, \textit{one more than} the greatest depth of an
essential $\mathcal{J}$-class,
Let $n$ denote the depth of the semigroup $S$. For each essential
$\mathcal{J}$-class $J$, let $G_J$ denote the maximal
subgroup of $J$. Now for every integer $0 \leq i < n$, let $K_i$ be the
direct product over all essential $\mathcal{J}$-classes of depth $i$ of $G_i$.

\begin{theorem} \textsf{(Depth Decomposition Theorem, Eilenberg-Tilson 1976)}\
Let $S$ be a finite semigroup of depth $n$, and let $K_0,
\dots, K_{n-1}$ be as defined above. Then there exist aperiodic
monoids $A_0, \dots A_n$ such that $S$ divides the wreath product
$$A_n \wr K_{n-1} \wr A_{n-1} \wr \dots \wr K_0 \wr A_0.$$
\end{theorem}
Thus, the depth decomposition theorem gives, for any finite
semigroup $S$, a Krohn-Rhodes decomposition with group length equal
to the depth of $S$. To apply the depth decomposition theorem, we need some information
about the $\mathcal{J}$-class structure and maximal subgroups of our
semigroups. The following proposition provides a description;
various parts of it have been observed
before \cite{Almeida05,Okninski98,Putcha88,Volkov03} but for
completeness we prove the entire statement.
\begin{proposition}\label{prop_tutdepth}
Let $n$ be a positive integer and $k$ a finite field. Then
\begin{itemize}
\item[(i)] $T_n(k)$ has
depth $n-1$ if $k = \BZ_2$, or depth $n$ otherwise. For $0 \leq i
\leq n-2$ or $0 \leq i \leq n-1$ as appropriate, $T_n(k)$ has
$\binom{n}{i}$ essential $\mathcal{J}$-classes of depth $i$, each of
which has maximal subgroup isomorphic to $T_{n-i}^*(k)$;
\item[(ii)] $UT_n(k)$ has depth $n-1$. For $0 \leq i < n-1$, $UT_n(k)$
has $\binom{n}{i}$ essential $\mathcal{J}$-classes of depth $i$,
each of which has maximal subgroup isomorphic to $UT_{n-1}^*(k)$;
and
\item[(iii)] $PT_n(k)$ has depth $n-1$. For $0 \leq i < n-1$, $PT_n(k)$ has
$\binom{n}{i}$ essential $\mathcal{J}$-classes of depth $i$, each of
which has maximal subgroup isomorphic to $PT_{n-i}^*(k)$.
\end{itemize}
\end{proposition}
\begin{proof}
We begin with the case of $T_n(k)$. By
Proposition~\ref{prop_regular_characterization}, the regular
$\mathcal{J}$-classes are exactly the $\mathcal{J}$-classes of the
subidentites. Moreover, if $e$ and $f$ are two subidentities, it is
easily seen (for example, by using
Proposition~\ref{prop_greens_characterization}), that $e$ is
$\mathcal{J}$-below $f$ if and only if $ef=fe =e$. Thus, the lattice
of regular $\mathcal{J}$-classes is isomorphic to the lattice
$\{0,1\}^n$, that is to the subset lattice of the set $\lbrace 1,
\dots, n \rbrace$. In particular, there are $\binom{n}{i}$ regular
$\mathcal{J}$-classes at depth $i$ for $i \in \lbrace 0, \dots, n
\rbrace$.

Now let $e \in T_n(k)$ be a subidentity at depth $i$, so that $e$ has
rank $n-i$.  It is easily seen that $eT_n(k)e$ is isomorphic to
$T_{n-i}(k)$ via the map that removes from a matrix all rows and
columns for which $e$ has a zero in the corresponding diagonal
position. Thus the maximal subgroup at $e$ is isomorphic to
$T_{n-i}^*(k)$.

Hence, in the case that $k \neq \BZ_2$, all regular $\mathcal{J}$-classes except
for that of $0$ are essential, giving the required result. In the case that
$k = \BZ_2$, however, $T_1^*(k)$ is trivial and so there are no essential
$\mathcal{J}$-classes of depth $n-1$. Thus, in this case, the depth of the
semigroup is one less.

The case of the unitriangular semigroup $UT_n(k)$ is exactly the same
except that the maximal subgroup of the $\mathcal{J}$-class of a subidentity
with $n-i$ diagonal entries is isomorphic to the unitriangular group
$UT_{n-i}^*(k)$. However, since $UT_1^*(k)$ is trivial regardless of the
field $k$, there are never essential $\mathcal{J}$-classes
of depth $n-1$, so the semigroup has depth $n-1$.

For the projective triangular semigroups $PT_n(k)$,
Proposition~\ref{prop_ptstructure} tells us that the lattice of
$\mathcal{J}$-classes is the same as that of $T_n(k)$; the maximal
subgroup of the $\mathcal{J}$-classes of a subidentity of rank $n-i$
is clearly the projective image $PT_{n-1}^*(k)$ of $T_{n-i}^*(k)$. In
particular, $PT_1^*(k)$ is trivial so as in the unitriangular case
there are
no essential $\mathcal{J}$-classes of depth $n-1$, and the semigroup
has depth $n-1$.
\end{proof}
Proposition~\ref{prop_tutdepth} supplies the information needed to apply
the Depth Decomposition Theorem to our semigroups. Doing so, we obtain:
\begin{align*}
T_n(k) &\prec A_n \wr T_1^*(k)^n \wr A_{n-1} \wr T_2^*(k)^{\binom{n}{2}} \wr \dots \wr T_n^*(k) \wr A_0 \\
UT_n(k) &\prec B_{n-1} \wr UT_2^*(k)^{\binom{n}{2}} \wr B_{n-2} \wr UT_3^*(k)^{\binom{n}{3}} \wr \dots \wr UT_n^*(k) \wr B_0 \\
PT_n(k) &\prec C_{n-1} \wr PT_2^*(k)^{\binom{n}{2}} \wr C_{n-2} \wr
PT_3^*(k)^{\binom{n}{3}} \wr \dots \wr PT_n^*(k) \wr C_0
\end{align*}
for some aperiodic semigroups $A_0, \dots, A_n, B_0, \dots, B_{n-1}, C_0, \dots, C_{n-1}$.
Thus, depth decomposition gives alternative (by Theorem~\ref{thm_trimatcomplexity}, optimal)
decompositions of group length $n-1$ for $UT_n(k)$ and $PT_n(k)$ and a
(suboptimal) group length $n$ decomposition for $T_n(k)$.
The theorem as stated does not give an explicit description of the
aperiodic terms; however, the interested reader could compute
appropriate ones through an analysis of the proof \cite{TilsonXI}.

\section*{Acknowledgements}

The research of the first author was supported by a Marie Curie
Intra-European Fellowship within the 6th European Community
Framework Programme. The first author would also like to thank
Kirsty for all her support and encouragement.  The work of the
second author was supported by an NSERC discovery grant.

\def\cprime{$'$}

\end{document}